\input amstex
\documentstyle{amsppt}

\define\Int{\operatorname{Int}}
\define\IntR1{\operatorname{R1-Int}}
\define\InvR1{\operatorname{R1-Inv}}
\define\I{\operatorname{I}}

\topmatter

\title
THE UNIQUENESS OF THE SOLUTION TO INVERSE INTERPOLATION PROBLEMS
\endtitle

\footnote""{This note is a revised and translated form of [6].}

\rightheadtext{INVERSE INTERPOLATION PROBLEMS}

\author
L.~V.~Veselova and O.~E.~Tikhonov
\endauthor

\address
Department of Higher Mathematics, Kazan State Technological 
University, Karl Marx Str.\ 68, Kazan, Tatarstan, 
420015 Russia
\endaddress

\address
Research Institute of Mathematics and Mechanics, 
Kazan State University, Universitetskaya Str.\ 17, Kazan, Tatarstan, 
420008 Russia
\endaddress

\email
Oleg.Tikhonov{\@}ksu.ru
\endemail

\keywords
Banach couple, interpolation space, exact interpolation space
\endkeywords 

\subjclass
Primary 46B70; Secondary 46M35
\endsubjclass

\thanks
Supported by the Russian Foundation for Basic 
Research, grants no.\ 95--01--00025 and no.\ 98--01--00103.
\endthanks

\abstract
We prove that an arbitrary Banach couple is uniquely determined by the 
collection of intermediate spaces that are interpolation spaces 
for the operators of rank one. From this we deduce a confirmation to the 
conjecture by Yu.~A.~Brudny{\u \i} and N.~Ya.~Kruglyak that 
a Banach couple is uniquely determined by the collection of all its 
interpolation spaces. Some relevant problems are also analyzed. 

\endabstract
\endtopmatter

\document

\head
Introduction
\endhead

The principal result of this paper is an exhaustive answer to a 
question which originates from the classical paper by N. Aronszajn 
and E. Gagliardo [1] (see also [3, \S 1, Subsection 9]): To what 
extent is a Banach couple determined by its interpolation spaces? 
Namely, we prove that an arbitrary Banach couple is uniquely 
determined by the collection of all its interpolation spaces. Note 
that Yu.\ A. Brudny{\u \i} and N. Ya.\ Kruglyak conjectured that this 
holds true and posed the corresponding problem [4, Conjecture 2.2.32 
and Problem 2.7.4(a)]. 

Also, we examine a similar problem for exact interpolation. The 
two problems are closely connected: the exhaustive solution to the 
second is based on that to the first. However, as what we mean, these 
are rather distinct. For example, in the case when the spaces of a 
Banach couple coincide as linear ones, the first problem is trivial in 
contrast to that solving the second has required certain efforts (see 
[5]).

We examine the first problem in Section 2 and the second one in 
Section 3. In each of the cases, we first analyze the corresponding 
problem for operators of rank one, and then we obtain the main result 
of the corresponding section as a corollary. 
The main results are already published in [7] and [8], but in contrast 
to those papers we underline the role of rank one operators here. 

\head
1. Notation and Preliminaries
\endhead

In this section, we introduce notation and recall some definitions 
and facts from the theory of interpolation of linear operators in 
Banach spaces (see, e.g., [2], [3], or [4]). Also, we present the 
proofs of some preliminary results.

For a Banach space $E$, we will denote by $\| \cdot \| _E$ the norm 
in $E$, by $B_E$ and by $B_E ^{\circ }$ the closed and the open unit 
ball in $E$ respectively, and by $\overline{D} ^E$ the closure in 
$E$ of a subset $D \subset E$. We will denote by $\| \cdot \| ^E$ 
the norm in the conjugate space $E^*$. 

Let $E$ and $F$ be two Banach spaces with $E \hookrightarrow F$. 
(Here and subsequently, the notation $E \hookrightarrow F$ for two 
Banach spaces $E$ and $F$ means that $E$ is embedded linearly and 
continuously into $F$.) The {\it embedding constant\/}, 
$\alpha (F,E)$, is defined by  
$$
\alpha (F,E) = 
\sup \bigl \{ \, \| x \| _F / \| x\| _E : 
x \in E \setminus \{ 0 \} \, \bigr \} .
$$
If $\varphi \in F^*$ then, clearly, $\varphi | _E \in E^*$. We write 
$\| \varphi \| ^E$ to denote $\| \varphi | _E \| ^E$ and put 
$$
\beta (E^*,F^*) = \sup \bigl \{ \, \| \varphi \| ^E / \| \varphi \| ^F
 : \varphi \in F^* \setminus \{ 0 \} \, \bigr \} . 
$$
The restriction operator $\varphi \mapsto \varphi | _E$ from 
$F^*$ into $E^*$ is adjoint to the embedding operator $E \to F$, 
therefore $\alpha (F,E) = \beta (E^*,F^*)$. 

\proclaim{Lemma 1.1}
Let $E$ and $F$ be Banach spaces with $E \hookrightarrow F$. 
Suppose that $E$ and $F$ do not coincide as linear spaces. Then 
there exists a sequence $\{ \varphi _n \} \subset F^*$ such that 
$\| \varphi _n \| ^F = 1$ and $\| \varphi _n \| ^E \to 0$. 
\endproclaim

\demo{Proof} 
It is well-known (e.g., [2, Lemma I.1.1]) 
that the unit ball $B_E$ is nowhere dense in $F$. 
Therefore $\overline{nB_E}^F$ does not contain $B_F$ for any natural 
number $n$. Hence there exists $x_n \in B_F$ such that 
$x_n \notin \overline{nB_E}^F$. It follows from the Hahn---Banach 
theorem that there exists $\varphi _n \in F^*$ such that  
$\| \varphi _n \| ^F = 1$ and 
$$
\sup \{ | \varphi _n(x) | : x \in \overline{nB_E} ^F \} < 
\varphi _n(x_n).
$$ 
Hence we have 
$$
\| \varphi _n \| ^E = \sup \{ | \varphi _n (x) | : x \in B_E \} 
\le \varphi _n(x_n) / n \le 1/n.
$$ 
\enddemo 

For two Banach spaces $E$ and $F$, the notation $E \simeq F$ will 
stand for the case when $E$ and $F$ coincide as linear spaces and 
their norms are equivalent, the notation $E\cong F$ will stand for 
the case when the norms are proportional, and $E=F$ will stand for 
the case when the norms coincide. 

Two Banach spaces $X$ and $Y$ are said to form a {\it Banach couple} 
$(X,Y)$ if they both are linearly and continuously embedded into a 
certain Hausdorff topological vector space. Note that if $X$ and $Y$ 
form a Banach couple and one of them is embedded into another as a 
set, then the embedding operator is linear and continuous [2, Lemma 
I.3.3], and consequently if $X$ and $Y$ coincide as sets, then 
$X \simeq Y$. A Banach couple $(X,Y)$ is called {\it embedded} if 
one of the spaces is embedded into another. 

To each Banach couple $(X,Y)$ can be canonically associated two 
Banach spaces, the {\it intersection } $X \cap Y$ and the {\it sum } 
$X+Y$, with norms defined by 
$$
\| z \| _{X \cap Y} = \max \{ \| z \| _X , \| z \| _Y \}
\qquad (z \in X \cap Y) 
$$
and 
$$
\| z \| _{X + Y} = \inf \{ \| x \| _X + \| y \| _Y :
 x \in X, \, y \in Y, \, z=x+y \}
\qquad (z \in X+Y).
$$ 
A Banach space $Z$ is said to be {\it intermediate } for a Banach 
couple $(X,Y)$ if 
$X \cap Y \hookrightarrow Z \hookrightarrow X+Y$. 
The set of all intermediate spaces for a Banach couple  $(X,Y)$ will 
be denoted by $\I(X,Y)$. 

We will use without reference the following assertions which hold 
true for any Banach couple $(X,Y)$:

a) $X \cap Y$ is dense in $X$ if and only if $Y$ is dense in $X+Y$; 

b) if $X \cap Y$ is dense in both $X$ and $Y$, then it is dense in 
$X+Y$, too.

\noindent [2; Ch.\ I, \S 3, p.\ 23, Corollary 1 and Corollary 2].

\medskip 
A Banach couple $(X,Y)$  is called {\it regular} if $X \cap Y$ is 
dense in both $X$ and $Y$. Clearly, $(X,Y)$ is regular if and only 
if $X \cap Y$ is dense in $X+Y$. 

For two Banach couples $(X,Y)$ and $(V,W)$, we will write  
$(X,Y) \sim (V,W)$ if $\I (X,Y) = \I (V,W)$. Obviously,  
$(X,Y) \sim (V,W)$ if and only if $X \cap Y \simeq V \cap W$ and  
$X+Y \simeq V+W$.

\proclaim{Lemma 1.2 {\rm (see [1, the proof of Corollary 10.XIV])}} 
Suppose that for two Banach couples, $(X,Y)$ and $(V,W)$, the following 
conditions are fulfilled: 

{\rm 1)} $(X,Y) \sim (V,W)$, 

{\rm 2)} $X \hookrightarrow V$, 

{\rm 3)} $Y \hookrightarrow W$. 

\noindent Then $X \simeq V$ and $Y \simeq W$.
\endproclaim

\proclaim{Lemma 1.3 {\rm (e.g., [2, Lemma I.3.4])}} 
Let, for a Banach couple $(X,Y)$, a subset $D \subset X \cap Y$ be 
dense in $X$ and an inequality $ \| x \| _Y \le c \| x \|_X$ hold 
on $D$. Then $X$ is embedded into $Y$ with an embedding constant less 
than or equal to $c$. 
\endproclaim

Next, consider certain families of linear operators in a Banach couple 
and relevant classes of the intermediate spaces.

A linear operator $T$ in $X+Y$ is called a {\it bounded linear 
operator on a Banach couple} $(X,Y)$ if $T$ maps $X$ into $X$ and 
$Y$ into $Y$ boundedly. Let $L((X,Y))$ stand for the set of all such 
operators. $L((X,Y))$ is a Banach space with respect to the norm 
$$
\| T \| _{L((X,Y))}= \max \{ \| T \| _{X \to X}, \| T \| _{Y \to Y} \} .
$$
Let $R1((X,Y))$ stand for the subset of the operators of rank one.

A Banach space $Z \in \I(X,Y)$ is called an {\it interpolation space} 
relative to a Banach couple $(X,Y)$ if $T(Z) \subset Z$ for every 
$T \in L((X,Y))$. The collection of all interpolation spaces for a 
Banach couple $(X,Y)$ will be denoted by $\Int(X,Y)$. It is 
well-known that if $Z \in \I (X,Y)$, $T \in L((X,Y))$, and 
$T(Z) \subset Z$, then $T |_Z$ is a bounded operator in $Z$. 
Moreover, if $Z \in \Int (X,Y)$, then there exists a constant 
$C>0$ such that the following ``interpolation inequality'' holds 
true for any $T \in L((X,Y))$: 
$$ 
\| T \| _{Z \to Z} \le C \| T \| _{L((X,Y))} .
$$ 
If the inequality holds with $C=1$, then $Z$ is said to be an 
{\it exact interpolation space}. The collection of all exact 
interpolation spaces for a Banach couple 
$(X,Y)$ will be denoted by $\Int _1 (X,Y)$. 

When one restricts himself to considering operators of rank one, he 
has to distinguish between the concept of invariant spaces and that of 
interpolation spaces. We will say that a space $Z \in \I (X,Y)$ is 
{\it invariant\/} for the operators of rank one on a Banach couple 
$(X,Y)$ if $T(Z) \subset Z$ for every $T \in R1((X,Y))$. The collection 
of all such spaces will be denoted by $\InvR1 (X,Y)$. If there 
exists, in addition, a constant $C>0$ such that  
$$ 
\| T \| _{Z \to Z} \le C \| T \| _{L((X,Y))}
$$ 
for any $T \in R1((X,Y))$, then we will say that $Z$ is an {\it 
interpolation space for the operators of rank one\/}. The 
collection of all such spaces will be denoted by $\IntR1 (X,Y)$. If 
the last inequality holds with $C=1$, then we will say that $Z$ 
is an {\it exact interpolation space for the operators of rank 
one,\/} and we will denote the collection of all such spaces by 
$\IntR1 _1 (X,Y)$. 

It is immediate from the above definitions that 
$$
\Int _1 (X,Y) \subset \Int (X,Y) \subset \I (X,Y),
$$
$$
\IntR1 _1 (X,Y) \subset \IntR1 (X,Y) \subset \InvR1 (X,Y) \subset 
\I (X,Y), 
$$
$$
\Int (X,Y) \subset \IntR1 (X,Y) ,
$$
and
$$
\Int _1 (X,Y) \subset \IntR1 _1 (X,Y).
$$
Note also that 
$ X \cap Y , \, X+Y \in \Int _1 (X,Y) $. 

\proclaim{Proposition 1.4} 
Let $(X,Y)$ and $(V,W)$ be two Banach couples. Then 
$$
{\goth I} (V,W) \subset {\goth I} (X,Y) \iff 
V,W \in {\goth I} (X,Y) ,
$$
where ${\goth I}$ may be substituted by any of the symbols $\I$, 
$\InvR1$, $\IntR1$, $\Int$, $\IntR1 _1$, or $\Int _1$. 
\endproclaim

\demo{Proof} 
Let us examine the case ${\goth I} =\InvR1$. The implication 
$\implies $ is trivial. Let us prove the $\impliedby $. Suppose that
$V,W \in \InvR1 (X,Y)$. Let $Z \in \InvR1 (V,W)$ and 
$T \in R1((X,Y))$. Then $T$ maps $V$ into $V$ and $W$ into $W$ 
boundedly. It follows that $T | _{V+W} \in R1((V,W))$, hence 
$T(Z) \subset Z$. Thus $Z \in \InvR1  (X,Y)$. 

Consider another case, ${\goth I} = \Int _1$ for instance.
Let us prove the $\impliedby $. Suppose that 
$V,W \in \Int _1 (X,Y)$. Let $Z \in \Int _1 (V,W)$ and 
$T \in L((X,Y))$. We have $\| T \| _{V \to V} \le \| T \| _{L((X,Y))}$ 
and $\| T \| _{W \to W} \le \| T \| _{L((X,Y))}$. It follows that 
$\| T \| _{Z \to Z} \le \| T \| _{L((V,W))} \le \| T \| _{L((X,Y))}$. 
Thus $Z \in \hbox{Int} _1 (X,Y)$.

The rest of the cases are treated similarly.
\enddemo

Let $(X,Y)$ and $(V,W)$ be two Banach couples. We write 
$(X,Y) \simeq (V,W)$ if either $X \simeq V$ and $Y \simeq W$ or 
$X \simeq W$ and $Y \simeq V$. 
We write $(X,Y) \cong (V,W)$ if either $X \cong V$ and 
$Y \cong W$ or $X \cong W$ and $Y \cong V$. 

Let ${\goth I}$ stand for one of the symbols: $\Int$, $\InvR1$,
$\IntR1 $, or $\I$. We say that a Banach couple $(X,Y)$ {\it is 
uniquely determined by a collection ${\goth I}(X,Y)$} if the 
implication 
$$
{\goth I} (X,Y) = {\goth I} (V,W) \implies (X,Y) \simeq (V,W)
\eqno{(1)}
$$ 
holds true for any Banach couple $(V,W)$. 
Let ${\goth U}({\goth I})$ stand for the class of all such 
couples $(X,Y)$.

Now, let ${\goth I}$ stand for either $\Int _1$ or $\IntR1 _1$.
We say that a Banach couple $(X,Y)$ {\it is uniquely determined by 
a collection ${\goth I}(X,Y)$} if the implication 
$$ 
{\goth I} (X,Y) = {\goth I} (V,W) \implies (X,Y) \cong (V,W)
\eqno{(2)}
$$
holds true for any Banach couple $(V,W)$.
Similarly to the above, we will denote by ${\goth U}({\goth I})$ 
the class of all such couples $(X,Y)$.

It is clear that the inverse implications in (1) and (2) are always 
hold true. 

\proclaim{Proposition 1.5}

{\rm a)} ${\goth U} ( \I) \subset {\goth U} ( \InvR1 ) \subset 
{\goth U} ( \IntR1 ) \subset {\goth U} ( \Int )$. 

{\rm b)} ${\goth U} ( \IntR1 _1 ) \subset {\goth U} ( \Int _1)$.
\endproclaim

\demo{Proof} 
Let us prove that 
${\goth U} ( \InvR1 ) \subset {\goth U} ( \IntR1 )$. Let
$(X,Y) \in {\goth U} ( \InvR1 )$ and let, for a Banach couple $(V,W)$, 
it hold
$$
\IntR1  (V,W) = \IntR1  (X,Y).
$$
By making use of Proposition 1.4, we get 
$$
\InvR1 (V,W) = \InvR1 (X,Y).
$$
Hence $(V,W) \simeq (X,Y)$. Thus $(X,Y) \in {\goth U} ( \IntR1 )$.

The other inclusions can be analyzed similarly.
\enddemo

\remark{Remark} 
A Banach couple $(X,Y)$ with $X \simeq Y$ gives a trivial example of 
a couple from ${\goth U} ( \I)$.
\endremark

\head 
2. The Uniqueness of the Solution to Inverse Problem \\ 
of Interpolation of Linear Operators
\endhead

Throughout this section, we will denote, for short, $X \cap Y$ and 
$X+Y$ by $\Delta $ and $\Sigma $, respectively. Also, we will denote 
by $\overline{Z}$ the closure in $\Sigma $ of a space 
$Z \in \I (X,Y)$. 

First, let us study when a Banach couple $(X,Y)$ is uniquely determined 
by the collection $\InvR1 (X,Y)$.

\proclaim{Proposition 2.1} 
Let $(X,Y)$ be a Banach couple and $Z \in \I (X,Y)$. 
Then $Z \in \InvR1 (X,Y)$ if and only if at least one of the 
following four conditions is satisfied: 

{\rm (i)} 
$\Delta \hookrightarrow Z \hookrightarrow \overline{\Delta}$,

{\rm (ii)} $X \hookrightarrow Z \hookrightarrow \overline{X}$,

{\rm (iii)} $ Y \hookrightarrow Z \hookrightarrow \overline{Y}$,

{\rm (iv)} $Z \simeq \Sigma$.
\endproclaim

\demo{Proof} 
By [1, Theorem 7.IV], if $Z \in \Int (X,Y)$ then at least one of the 
four conditions is satisfied. However, it was actually proved there 
that one of (i)--(iv) was satisfied whenever $Z \in \InvR1 (X,Y)$. 
Therefore, it remains to show that each of the conditions entails that 
$Z \in \InvR1 (X,Y)$. 

Let $T$ be an operator of rank one in $\Sigma$, which maps boundedly 
$X$ into $X$ and $Y$ into $Y$. Necessary, $T$ is of the form 
$T( \cdot )= \varphi ( \cdot ) x$ with $\varphi \in {\Sigma }^*$, 
$x \in \Sigma$. The following two cases are possible: $x \in \Delta$ 
and $x \notin \Delta$. In the first case, $T$ maps boundedly $Z$ 
into $Z$ for every $Z \in \I (X,Y)$. Now, let $x \notin \Delta$. We 
will show that if one of the conditions (i)--(iv) is fulfilled 
then $T$ maps $Z$ into itself. 

First, let (i) be satisfied. Since $T$ has to map $\Delta$ into 
itself, $\varphi | _{\Delta} =0$, therefore 
$\varphi | _Z = 0$, and hence $T$ maps $Z$ into itself.

Next, let (ii) be satisfied. If $x \in X$ then, clearly, $T$ maps $Z$ 
into itself. If $x \notin X$ then $\varphi | _X = 0$, therefore 
$\varphi | _Z =0$, and hence $T$ maps $Z$ into itself. 

The case (iii) is similar to (ii), and (iv) is trivial.
\enddemo

\remark{Remark} 
The proposition shows that a Banach couple is not determined 
in general by the collection of all its intermediate spaces which 
are invariant for the operators of rank one. In fact, if $(X,Y)$ is 
a regular Banach couple, then 
$ \InvR1 (X,Y) = \InvR1  (\Delta, \Sigma) = \I (X,Y)$, and if the 
couple is non-embedded, then $(X,Y) \not \simeq (\Delta, \Sigma)$.
\endremark

\proclaim{Lemma 2.2} 
Suppose that for two Banach couples, $(X,Y)$ and $(V,W)$, the 
following four conditions are fulfilled:

{\rm 1)} $X \hookrightarrow Y$, 

{\rm 2)} $X$ is not dense in $Y$, 

{\rm 3)} $(X,Y) \sim (V,W)$, 

{\rm 4)} $ \InvR1 (V,W) \subset \InvR1  (X,Y)$. 

\noindent Then $(X,Y) \simeq (V,W)$.
\endproclaim 

\demo{Proof} 
To prove the lemma, it suffices to examine all the positions of 
$V,W \in \InvR1 (X,Y)$ relative to $(X,Y)$, which are admissible by 
Proposition 2.1.
\enddemo

\proclaim{Lemma 2.3}
Suppose that the following three relations hold for a Banach couple 
$(X,Y)$ and a non-embedded Banach couple $(V,W)$: 

{\rm 1)} $\Delta$ is dense in neither $X$ nor $Y$, 

{\rm 2)} $(X,Y) \sim (V,W)$, 

{\rm 3)} $ \InvR1 (V,W) \subset \InvR1 (X,Y)$. 

\noindent Then $(X,Y) \simeq (V,W)$. 
\endproclaim 

\demo{Proof}
First, note that neither $X$ nor $Y$ is dense in $\Sigma$. 
If one examines all admissible positions of $V,W \in \InvR1 (X,Y)$ 
relative to $(X,Y)$ and rejects the cases which entail the embeddedness 
of $(V,W)$, then he easily sees that it suffices to 
consider the four cases only:

(i) $\Delta \hookrightarrow V \hookrightarrow \overline{\Delta}$, \,\,\,
    $\Delta \hookrightarrow W \hookrightarrow \overline{\Delta}$;

(ii) $X \hookrightarrow V \hookrightarrow \overline{X}$, \,\,\,
     $X \hookrightarrow W \hookrightarrow \overline{X}$;

(iii) $\Delta \hookrightarrow V \hookrightarrow \overline{\Delta}$, \,\,\,
      $X \hookrightarrow W \hookrightarrow \overline{X}$;

(iv) $X \hookrightarrow V \hookrightarrow \overline{X}$, \,\,\,
     $Y \hookrightarrow W \hookrightarrow \overline{Y}$.

\noindent Each of the cases (i), (ii), and (iii) leads to 
$V+W \hookrightarrow \overline{X}$, but this contradicts the fact 
that $V+W \simeq \Sigma$ and $\overline{X} \not \simeq \Sigma$. 
By Lemma 1.2, (iv) implies  $(X,Y) \simeq (V,W)$, which completes 
the proof.
\enddemo

It is easy to deduce the following theorem from Lemma 2.2 and Lemma 
2.3.

\proclaim{Theorem 2.4} 
Suppose that a non-regular Banach couple satisfies one of the two 
conditions:

{\rm a)} $(X,Y)$ is embedded,

{\rm b)} $\Delta$ is dense in neither $X$ nor $Y$.

\noindent Then $(X,Y) \in {\goth U} (\InvR1 )$. 
\endproclaim

We now turn to proving the fact that an arbitrary Banach couple $(X,Y)$ 
is uniquely determined by the collection $\IntR1 (X,Y)$.

\proclaim{Lemma 2.5} 
Suppose that for two Banach couples, 
$(X,Y)$ and $(V,W)$, the following conditions are satisfied: 

{\rm 1)} $(X,Y) \sim (V,W)$,

{\rm 2)} $\Delta$ is dense in $Y$, 

{\rm 3)} $X \hookrightarrow V$,

{\rm 4)} $V \not \simeq \Sigma$,

{\rm 5)} $W \in \IntR1 (X,Y)$.

\noindent Then $X\simeq V$ and $Y \simeq W$.
\endproclaim

\demo{Proof} 
By Lemma 1.1, we can choose a sequence 
$\{ \varphi _n \} \subset \Sigma ^*$ such that 
$\| \varphi _n \| ^{\Sigma} = 1$ and $\| \varphi _n \| ^V \to 0$. 
Then from the condition 3) we obtain $\| \varphi _n \| ^X \to 0$.
By making use of the well-known relation 
$B_{\Sigma} ^{\circ} = \hbox{conv}(B_X ^{\circ}, B_Y ^{\circ} )$, 
we conclude that $\| \varphi _n \| ^Y \to 1$. 
Since $\Sigma \simeq V+W$, there exists a constant $c>0$ such that 
$\| \varphi _n \| ^{V+W} \ge c$ for all $n$. 
Since $\| \varphi _n \| ^V \to 0$, it follows that there exist a 
constant $c_1 >0$ and a natural $n_0$ such that
$\| \varphi _n \| ^W \ge c_1$ as $n \ge n_0$. 

By 5), there exists a constant $C>0$ such that
$$ 
\| T \| _{W \to W} \le C \| T \| _{L((X,Y))}
$$ 
for any $T \in R1((X,Y))$. Let us take an arbitrary $x \in \Delta$ 
and consider the linear operators $T_n$ of rank one given by the 
formula
$T_n( \cdot )=\varphi _n( \cdot )x$. Then we have
$$
\| T_n \| _{X \to X} = \| \varphi _n \| ^X \| x \| _X \to 0,
$$
$$
\| T_n \| _{Y \to Y} = \| \varphi _n \| ^Y \| x \| _Y \to \| x \| _Y,
$$
$$
\| T_n \| _{W \to W} = \| \varphi _n \| ^W \| x \| _W 
\ge c_1 \| x \| _W \quad \text{as} \quad n \ge n_0 ,
$$
It follows that if we take $n$ sufficiently large then we have 
$$
c_1 \| x \| _W \le \| T_n \| _{W \to W} \le C 
\max \{ \| T_n \| _{X \to X}, \| T_n \| _{Y \to Y} \} \le 
2C \| x \| _Y .
$$
Hence $\| x \| _W \le (2C/c_1) \| x \| _Y$ for all $x \in \Delta$. 
Since $\Delta$ is dense in $Y$, we have $Y \subset W$ by Lemma 1.3. 
Now, application of Lemma 1.2 completes the proof.
\enddemo

\proclaim{Theorem 2.6} 
Suppose that for two Banach couples, 
$(X,Y)$ and $(V,W)$, the following conditions are satisfied: 

{\rm 1)} $(X,Y) \sim (V,W)$, 

{\rm 2)} each of the couples $(X,Y)$ and $(V,W)$ is non-embedded, 

{\rm 3)} $\IntR1 (V,W) \subset \IntR1 (X,Y)$.

\noindent Then $(X,Y) \simeq (V,W)$.
\endproclaim

\demo{Proof} 
To prove the theorem, we analyze the three cases separately.

a) 
Suppose, first, that $\Delta$ is dense in neither $X$ nor $Y$. In 
this case, the assertion of the theorem follows from Lemma 2.3, 
since $\IntR1 (V,W) \subset \IntR1 (X,Y)$ implies 
$\InvR1 (V,W) \subset \InvR1 (X,Y)$ (see Proposition 1.4).

b) 
Next, suppose that $\Delta$ is dense in $Y$, but not in $X$. Recall 
that in this case $X$ is dense in $\Sigma $, but $Y$ is not. We 
examine all the positions of $V,W \in \InvR1 (X,Y)$ relative to 
$(X,Y)$, which are admissible by Proposition 2.1, and reject the 
positions which entail that $(V,W)$ is embedded. We see that it 
suffices to treaty the six cases: 

(b1) $\Delta \hookrightarrow V \hookrightarrow \overline{\Delta}$,
\,\,\, $\Delta \hookrightarrow W \hookrightarrow \overline{\Delta}$;

(b2) $Y \hookrightarrow V \hookrightarrow \overline{Y}$, \,\,\,
   $\Delta \hookrightarrow W \hookrightarrow \overline{\Delta}$;

(b3) $X \hookrightarrow V \hookrightarrow \overline{X}$, \,\,\,
   $\Delta \hookrightarrow W \hookrightarrow \overline{\Delta}$;

(b4) $Y \hookrightarrow V \hookrightarrow \overline{Y}$, \,\,\,
   $Y \hookrightarrow W \hookrightarrow \overline{Y}$;

(b5) $X \hookrightarrow V \hookrightarrow \overline{X}$, \,\,\,
   $X \hookrightarrow W \hookrightarrow \overline{X}$;

(b6) $X \hookrightarrow V \hookrightarrow \overline{X}$, \,\,\,
   $Y \hookrightarrow W \hookrightarrow \overline{Y}$.

\noindent Each of the cases (b1), (b2), and (b4) entails that 
$V+W \hookrightarrow \overline{Y}$, but this contradicts the fact 
that $V+W \simeq \Sigma $ and $\overline{Y} \not \simeq \Sigma $. 
The case (b5) yields $X \hookrightarrow V \cap W$, which 
contradicts the fact that $V \cap W \simeq \Delta $ and 
$\overline{\Delta } ^X \not \simeq X$. By Lemma 1.2, (b6) gives 
$X \simeq V$ and $Y \simeq W$. It remains to treaty (b3). In this 
case, all the assumptions of Lemma 2.5 turn out to be satisfied. 
Therefore, $(X,Y) \simeq (V,W)$. 

c) 
Suppose, finally, that $\Delta$ is dense in both $X$ and $Y$. 
We will analyze two possibilities, c1) and c2).

c1) 
Let $Y+V \simeq \Sigma $. We proceed similarly to the proof of Lemma 
2.5. We choose a sequence $\{ \varphi _n \} \subset \Sigma ^*$ such 
that $\| \varphi _n \| ^{\Sigma} = 1$ and 
$\| \varphi _n \| ^Y \to 0$. Then $\| \varphi _n \| ^X \to 1$ and 
there exist a constant $c>0$ and a natural $n_0$ such that $\| 
\varphi _n \| ^V \ge c$ as $n \ge n_0$. For $x \in \Delta $, we 
estimate the norms of the operators $T_n$ defined by 
$T_n ( \cdot ) = \varphi _n( \cdot )x$, in each of the spaces $V$, 
$X$, and $Y$. Then from $V \in \Int _1 (X, Y)$ we infer that there 
exists a constant  $c_1$ such that  $\| x \| _V \le c_1 \| x \| _X$ 
for all $x \in \Delta $. But $\Delta $ is dense in $X$, and 
consequently $X \subset V$. Now, as in the case (b3), all the 
assumptions of Lemma 2.5 turn out to be satisfied. Therefore, 
$(X,Y) \simeq (V,W)$. 

c2) 
Let $Y+V \not \simeq \Sigma $. Again, we proceed similarly to the 
proof of Lemma 2.5. Here we take a sequence 
$\{ \varphi _n \} \subset \Sigma ^*$ such that 
$\| \varphi _n \| ^{\Sigma} = 1$ and $\| \varphi _n \| ^{Y+V} \to 0$.
Then $\| \varphi _n \| ^Y \to 0$ and $\| \varphi _n \| ^V \to 0$. It 
follows that $\| \varphi _n \| ^X \to 1$ and that there exist a 
constant $c>0$ and a natural $n_0$ such that 
$\| \varphi _n \| ^W \ge c$ as $n \ge n_0$. For $x \in \Delta $, we 
estimate the norms of the operators 
$T_n( \cdot ) = \varphi _n( \cdot )x$ in $W$, $X$, and $Y$ 
and conclude that there exists a constant $c_1$ such that 
$\| x \| _W \le c_1 \| x \| _X$ for all $x \in \Delta $. From this 
we obtain $X \subset W$, and application of Lemma 2.5 (with 
interchanged $V$ and $W$) completes the proof.
\enddemo

\remark{Remark}
N.~Aronszajn and E.~Gagliardo said in [1, Remark 10.XV] that they 
did not know of any example of two Banach couples $(X,Y)$ and $(V,W)$ 
with $(X,Y) \not \simeq (V,W)$ and satisfying the following 
conditions:

1) $(X,Y) \sim (V,W)$, 

2) $(X,Y) \not \simeq (\Delta, \Sigma)$, 

3) $(V,W) \not \simeq (\Delta, \Sigma)$, 

4) $\Int (V,W) \subset \Int (X,Y)$.

\noindent Taking into account Proposition 1.4, one sees that if two 
Banach couples satisfy all of 1)--4) then they satisfy the 
assumptions of Theorem 2.6. Therefore, the realization of the 
conditions 1)--4) does yield $(X,Y) \simeq (V,W)$.
\endremark

\proclaim{Theorem 2.7} 
Let $(X,Y)$ be a regular non-embedded Banach couple. Then neither 
$X$ nor $Y$ belongs to $\IntR1 (\Delta, \, \Sigma)$. 
\endproclaim

\demo{Proof} 
Suppose on the contrary that $Y \in \IntR1 (\Delta , \Sigma )$, for 
instance. 

Again, we proceed similarly to the proof of Lemma 2.5.
We choose a sequence 
$\{ \varphi _n \} \subset \Sigma ^*$ such that 
$\| \varphi _n \| ^{\Sigma} = 1$ and $\| \varphi _n \| ^X \to 0$.
Then $\| \varphi _n \| ^Y \to 1$ and 
$\| \varphi _n \| ^{\Delta} \to 0$.  
For $x \in \Delta $, we estimate the norms of the operators 
$T_n(  \cdot )=\varphi _n( \cdot )x$ in each of the spaces $Y$, 
$\Delta $, and $\Sigma $. Taking into account the fact that 
$Y \in \IntR1 (\Delta , \Sigma )$, we infer that there exists a 
constant $c$ such that $\| x \| _Y \le c \| x \| _{\Sigma } $ for 
all $x \in \Delta $. But $\Delta $ is dense in $\Sigma $, and 
consequently $\Sigma \subset Y$. This contradicts the assumption 
of that $(X,Y)$ is not embedded.
\enddemo
 
By combining Theorem 2.4, Theorem 2.6, Theorem 2.7, and Proposition 
1.5, we obtain the following theorem and corollary, which concludes 
the section. 

\proclaim{Theorem 2.8}
Let $(X,Y)$ and $(V,W)$ be two Banach couples with 
$$
\IntR1 (X,Y) = \IntR1 (V,W).
$$
Then $(X,Y) \simeq (V,W)$.
\endproclaim

\proclaim{Corollary 2.9 {\rm [7]}}
Let $(X,Y)$ and $(V,W)$ be two Banach couples with 
$$
\Int (X,Y) = \Int (V,W).
$$
Then $(X,Y) \simeq (V,W)$.
\endproclaim

\head 
3. The Uniqueness of the Solution to Inverse Problem \\ 
of Exact Interpolation
\endhead

As proportional norms on a linear space define the same operator 
norm, we will replace norms by proportional ones when this is 
convenient. This will involve no loss of generality. In particular, 
embedding constants will be frequently assumed to equal 1.

\proclaim{Lemma 3.1} 
Suppose that two Banach couples, $(X,Y)$ and $(V,W)$, satisfy the 
following conditions:

{\rm 1)} $X \simeq V$, $Y \simeq W$;

{\rm 2)} $X \not \simeq X+Y$;

{\rm 3)} $W \in \IntR1 _1 (X,Y)$, $Y \in \IntR1_1 (V,W)$.

\noindent Then $Y \cong W$. 
\endproclaim

\demo{Proof} 
We examine two cases separately:
\par
a) $X$ is not dense in $X+Y$,
\par
b) $X$ is dense in $X+Y$.

a) 
Since $X$ is not dense in $X+Y$, we can find 
$ \varphi _0 \in (X+Y)^*$ such that  
$\| \varphi _0 \| ^{X+Y} = 1$ and $ \varphi _0 | _X = 0$. 
It is easy to deduce from 
$B_{X+Y}^{\circ} = \hbox{conv}(B_X^{\circ}, B_Y^{\circ})$   
that $\| \varphi _0 \| ^Y = 1$. 
Replacing the norm in $W$ by a proportional one if necessary, we may 
and do assume that $\| y_0 \| _Y = \| y_0 \| _W = 1$ for some
$y_0 \in Y$ ($\simeq W$). Consider the linear operator of rank one given 
by the formula $T_0( \cdot ) = \varphi _0( \cdot )y_0$. Then we have 
$$
\| T_0 \| _{X \to X} = \| T_0 \| _{V \to V} = 0,
$$
$$
\| T_0 \| _{Y\to Y} = \| \varphi _0 \| ^Y \| y_0 \| _Y = 1,
$$
$$
\| T_0 \| _{W\to W} = \| \varphi _0 \| ^W \| y_0 \| _W = 
\| \varphi _0 \| ^W. 
$$
By 3), it follows that $\| \varphi _0 \| ^W = 1$.
\par
For an arbitrary $y \in Y \simeq W$, consider the operator given by  
$S_0( \cdot ) = \varphi _0( \cdot )y$. We have 
$$
\| S_0 \| _{X \to X} = \| S_0 \| _{V \to V} = 0,
$$
$$
\| S_0 \| _{Y\to Y} = \| y \| _Y,
$$
$$
\| S_0 \| _{W \to W} = \| y \| _W .
$$
By 3), it follows that $\| y \| _Y = \| y \| _W$, i.e., $Y=W$.

b) 
Arguments similar to those in the proof of Lemma 2.5 show that we 
can choose a sequence $\{ \varphi _n \} \subset (X+Y)^*$ such that 
$\| \varphi _n \| ^X \to 0$ and $\| \varphi _n \| ^Y = 1$. 
Then, from 1), we obtain $\| \varphi _n \| ^V \to 0$ and 
$c_1 \le \| \varphi _n \| ^W \le c_2$ for some $c_1,c_2 > 0$.

Since $X$ is dense in $X+Y$, $X \cap Y$ is dense in $Y$. In 
particular, $X \cap Y \ne \{ 0 \}$. 
Without loss of generality we assume that 
$\| x_0 \| _Y = \| x_0 \| _W = 1$ for some 
$x_0 \in {X \cap Y}$ $( \simeq V \cap W)$. Then the norms of 
the rank one operators $T_n( \cdot ) = \varphi _n( \cdot )x_0$ 
satisfy
$$
\| T_n \| _{X \to X} = \| \varphi _n \| ^X \| x_0 \| _X \to 0,
$$
$$
\| T_n \| _{Y \to Y} = \| \varphi _n \| ^Y \| x_0 \| _Y = 1,
$$
$$
\| T_n \| _{V \to V} = \| \varphi _n \| ^V \| x_0 \| _V \to 0,
$$
$$
\| T_n \| _{W \to W} = \| \varphi _n \| ^W \| x_0 \| _W =
\| \varphi _n \| ^W. 
$$
By 3), it follows that $\| \varphi _n \| ^W = 1$ for all $n$ is large 
enough. 
\par
For an arbitrary $x \in X \cap Y$ and sufficiently large $n$'s, 
the norms of operators 
$S_n( \cdot ) = \varphi _n( \cdot )x$ satisfy 
$$
\| S_n \| _{X \to X} = \| \varphi _n \| ^X \| x \| _X \to 0,
$$
$$
\| S_n \| _{Y \to Y} = \| x \| _Y,
$$
$$
\| S_n \| _{V \to V} = \| \varphi _n \| ^V \| x \| _V \to 0,
$$
$$
\| S_n \| _{W \to W} = \| x \| _W.
$$
By 3), it follows that $\| x \| _Y = \| x \| _W$. Since $X \cap Y$ 
is dense in $Y$ ($\simeq W$), we have $Y=W$. 
\enddemo

\proclaim{Proposition 3.2} 
Suppose that an embedded Banach couple $(X,Y)$ and another Banach 
couple $(V,W)$ satisfy the following conditions: 

{\rm 1)} $(X,Y) \simeq (V,W)$,

{\rm 2)} $ \IntR1 _1 (X,Y) = \IntR1 _1 (V,W) $. 

\noindent Then $(X,Y) \cong (V,W)$. 
\endproclaim

The proof of the proposition is based on a sequence of lemmas which 
we present below. Throughout the lemmas and the proof of the 
proposition we assume that $X \hookrightarrow Y$ with 
$\alpha (Y,X) = 1$. Moreover, we fix a sequence $\{ u_n \} \subset X$ 
with $\| u_n \| _Y = 1$ and $\| u_n \| ^{-1}_X \to \alpha (Y,X) = 1$ 
as $n \to \infty $. We also fix a sequence 
$\{ \theta _n \} \subset Y^* $ with $\| \theta _n \| ^Y = 1$ and 
$\theta _n (u_n) = \| u_n \| _Y = 1$. By Lemma 3.3 below, 
$\| \theta _n \| ^X \to \beta (X^*,Y^*)=1$. 

\proclaim{Lemma 3.3}
Let $A$ be a bounded linear operator from a Banach space $E$ into a 
Banach space $F$. Suppose that we take $x_n \in E \setminus \{ 0 \}$ 
and $\varphi _n \in F^* \setminus \{ 0 \}$ such that 
$\| Ax_n \| _F / \| x_n \| _E \to \| A\| $ $(n \to \infty)$ and 
$\varphi _n (Ax_n)=\| \varphi _n \| ^F \| Ax_n \| _F$. Then 
$\| A^* \varphi _n \| ^E / \| \varphi _n \| ^F \to \| A^* \| $ 
$(n \to \infty)$. 
\endproclaim

\demo{Proof}
The proof is straightforward.
\enddemo

\proclaim{Lemma 3.4} 
Let $Z \in \IntR1 _1 (X,Y)$. Then the following relations hold:

{\rm (i)} 
$\| u_n \| _Y / \| u_n \| _Z \,\, (= \| u_n \| ^{-1}_Z) 
\,\to \alpha (Y,Z) $, 

{\rm (ii)} 
$\| \theta _n \| ^Z / \| \theta _n \| ^Y \,\, (= \| \theta _n \| ^Z) 
\, \to \beta (Z^*, Y^*)$, 

{\rm (iii)} 
$\| u_n \| _Z / \| u_n \| _X \to \alpha (Z,X)$,

{\rm (iv)} 
$\| \theta _n \| ^X / \| \theta _n \| ^Z \to \beta (X^*, Z^*)$.
\endproclaim

\demo{Proof} 
For an arbitrary $\varphi \in Y^* \setminus \{ 0 \}$, consider the 
linear operators $T_n$ of rank one defined by the formula 
$T_n ( \cdot ) = \varphi ( \cdot ) u_n$. Then we have 
$$
\| T_n \| _{Y \to Y} = \| \varphi \| ^Y \| u_n \| _Y = 
\| \varphi \| ^Y,
$$
$$
\| T_n \| _{X \to X} = \| \varphi \| ^X \| u_n \| _X \le 
\| \varphi \| ^Y \| u_n \| _X, 
$$
$$
\| T_n \| _{Z \to Z} = \| \varphi \| ^Z \| u_n \| _Z.
$$
Since $Z\in \IntR1 _1 (X,Y)$, it follows that
$$
\| \varphi \| ^Z \| u_n \| _Z \le \| \varphi \| ^Y \| u_n \| _X.
$$
Thus 
$\| u_n \| _X / \| u_n \| _Z \ge \| \varphi \| ^Z / \| \varphi \| ^Y$ 
for every $\varphi \in Y^* \setminus \{ 0 \}$, therefore 
$$
\| u_n \| _X / \| u_n \| _Z \ge \beta (Z^*,Y^*) = \alpha (Y,Z).
$$
Hence we obtain
$$
\alpha (Y,Z) \ge \frac{\| u_n \| _Y}{\| u_n \| _Z} = 
\frac{1}{\| u_n \| _Z} = 
\| u_n \| ^{-1}_X \frac{ \| u_n \| _X}{ \| u_n \| _Z} 
\ge \| u_n \| ^{-1}_X \alpha (Y,Z).
$$
Since $\| u_n \| ^{-1}_X \to 1$, it follows that 
$\| u_n \| ^{-1}_Z \to \alpha (Y,Z)$, and (i) is proved.
\par
The relation (ii) follows from (i) and Lemma 3.3.
\par 
The proof of (iv) is a slight modification of that of (i). Namely, 
instead of $T_n$ we consider for an arbitrary 
$x \in X \setminus \{ 0 \}$ the operators $S_n$ defined by
$S_n( \cdot ) = \theta _n( \cdot )x$. 
\par
To prove (iii), we note that 
$$
\alpha (Z,X) \alpha (Y,Z) = \beta (X^*,Z^*) \beta (Z^*,Y^*)
$$
$$
= \lim \frac {\| \theta _n \| ^X} {\| \theta _n \| ^Z} 
\lim \frac {\| \theta _n \| ^Z} {\| \theta _n \| ^Y} =
\lim \frac {\| \theta _n \| ^X} {\| \theta _n \| ^Y} = 1.
$$
Hence 
$$
\frac {\| u_n \| _Z } {\| u_n \| _X} =
\frac {\| u_n \| _Y / \| u_n \| _X} {\| u_n \| _Y / \| u_n \| _Z} 
\to \frac {1} {\alpha (Y,Z)} = \alpha (Z,X). 
$$
\enddemo

\proclaim{Lemma 3.5} 
Let $\IntR1 _1 (X,Y) = \IntR1 _1 (Z,Y)$. Then $X\cong Z$.
\endproclaim

\demo{Proof}
We can and do assume that $\alpha (Y,Z) = 1$. Clearly $X\simeq Z$. 
Applying part (iii) and then part (i)
of Lemma 3.4, we obtain for an arbitrary 
$x\in X \setminus \{ 0 \}$: 
$$
\frac {\| x \| _Z } {\| x \| _X} \le \alpha (Z,X) = 
\lim \frac {\| u_n \| _Z } {\| u_n \| _X} =
\lim \frac {\| u_n \| _Y / \| u_n \| _X} 
{\| u_n \| _Y / \| u_n \| _Z} = 1 
$$
Similarly $\| x \| _X / \| x \| _Z \le 1$. Thus $X=Z$.
\enddemo

\demo{Proof of Proposition 3.2} 
If $X$ and $Y$ do not coincide as linear spaces then it suffices 
to apply Lemma 3.1 and Lemma 3.5. 
\par
Now, let $X\simeq Y$. Then $X \simeq Y \simeq V \simeq W$ and also 
$X^* \simeq Y^* \simeq V^* \simeq W^*$. 
From (i) in Lemma 3.4, we have 
$\| u_n \| ^{-1}_V \to \alpha (Y,V)$ and 
$\| u_n \| ^{-1}_W \to \alpha (Y,W)$. 
By interchanging the roles of $X$ and $Y$ in Lemma 3.4, we see that 
there exists a sequence $\{ \eta _n \} \subset X^*$ such that
$$
\| \eta _n \| ^Y = 1, \qquad 
\frac{1}{ \| \eta _n \| ^X } \to \beta (Y^*,X^*), 
$$
$$
\frac{1}{\| \eta _n \| ^V } \to \beta (Y^*,V^*), \quad \text{ and }
\quad \frac{1}{ \| \eta _n \| ^W } \to \beta (Y^*,W^*).
$$ 
Consider the operators $T_n$ of rank one defined by the formula
$T_n( \cdot ) = \eta _n( \cdot )u_n$. It is easy to calculate their 
norms in each of the spaces $X$, $Y$, $V$, and $W$. Then, from 
$ \Int _1^N (X,Y) = \Int _1^N (V,W) $, we obtain
$$
\max \{ 1, \, \| u_n \| _X \| \eta _n \| ^X \} = 
\max \{ \| u_n \| _W \| \eta _n \| ^W, 
\, \| u_n \| _V \| \eta _n \| ^V \}. 
$$

Note that
$$
\lim \max \{ 1, \, \| u_n \| _X \| \eta _n \| ^X \} = 1, 
$$
for otherwise we can find $\varepsilon > 0$ and an infinite family of
indices $\{ n_k \} $ such that 
$\| u_{n_k} \| _X \| \eta _{n_k} \| ^X  > 1+ \varepsilon $. Hence
$$
\alpha (Y,X) \alpha (X,Y) = \alpha (Y,X) \beta (Y^*,X^*) = 
\lim ( \| u_n \| _X \| \eta _n \| ^X) ^{-1} < 1, 
$$
which is impossible.

Clearly $V$ or $W$ (say, $W$) satisfies the following:
There exists an increasing infinite sequence of indices 
$\{ n_k \} $ such that 
$$
\max \{ 1, \, \| u_{n_k} \| _X \| \eta _{n_k} \| ^X \} = 
\| u_{n_k} \| _W \| \eta _{n_k} \| ^W.
$$

We can now write 
$$
1 = \lim_{k \to \infty} \| u_{n_k} \| _W \| \eta _{n_k} \| ^W = 
1/( \alpha (Y,W) \beta (Y^*,W^*)) = 
1/( \alpha (Y,W) \alpha (W,Y)). 
$$
Thus $ \alpha (Y,W) \alpha (W,Y) = 1$, and hence $Y\cong W$. 
An application of Lemma 3.5 completes the proof.
\enddemo

\proclaim{Proposition 3.6} 
Suppose that two Banach couples, $(X,Y)$ and $(V,W)$, satisfy the 
conditions:

{\rm 1)} $(X,Y) \simeq (V,W)$, 
 
{\rm 2)} $\IntR1 _1 (X,Y) = \IntR1 _1 (V,W)$.

\noindent Then $(X,Y) \cong (V,W)$. 
\endproclaim 

\demo{Proof} 
The case of an embedded couple $(X,Y)$ has been analyzed in 
Proposition 3.2. In the case where $(X,Y)$ is not embedded, we have
$X \not \simeq X+Y$ and $Y \not \simeq X+Y$, therefore it suffices 
to apply Lemma 3.1 twice.
\enddemo

Finally, we present the complete analogues of Theorem 2.8 and 
Corollary 2.9.

\proclaim{Theorem 3.7}
Let $(X,Y)$ and $(V,W)$ be two Banach couples with 
$$
\IntR1 _1 (X,Y) = \IntR1 _1 (V,W). 
$$
Then $(X,Y) \cong (V,W)$.
\endproclaim

\demo{Proof} \hfill
Proposition 1.4 shows that if \hfill 
$\IntR1 _1 (V,W) = \IntR1 _1 (X,Y)$ \hfill 
then $\IntR1 (V,W) = \IntR1 (X,Y)$. By Theorem 2.8, the latter 
yields $(V,W) \simeq (X,Y)$, and Proposition 3.6 gives 
$(V,W) \cong (X,Y)$.
\enddemo

\proclaim{Corollary 3.8 {\rm [8]}} 
Let $(X,Y)$ and $(V,W)$ be two Banach couples with 
$$
\Int _1 (X,Y) = \Int _1 (V,W). 
$$
Then $(X,Y) \cong (V,W)$.
\endproclaim

\demo{Proof}
It suffices to apply Proposition 1.5 b).
\enddemo

\remark{Concluding Remarks and Unsolved Problems}

1. The reformulated Theorem 2.8, Corollary 2.9, Theorem 3.7, and
Corollary 3.8 can be summarized as follows: Each of the 
${\goth U} (\Int _1)$, ${\goth U} (\Int )$, ${\goth U} (\Int _1 ^N)$, 
and ${\goth U} (\Int ^N)$ exhausts all Banach couples. 

2. In the remark that concludes Section 1, we have adduced the example of
a Banach couple from ${\goth U} (\I)$. Obviously, any non-embedded
Banach couple does not belong to ${\goth U} (\I )$. We believe that
to characterize the Banach couples from ${\goth U} (\I)$ is an
interesting and nontrivial problem.

3. Another problem that remains unsolved is to describe
${\goth U} (\InvR1 )$ (cf.\ Theorem 2.4 and Remark after
Proposition 2.1).

\endremark

\head
Acknowledgment
\endhead

We are indebted to N. M. Zobin who attracted our attention to inverse 
interpolation problems. We wish to express gratitude to Yu.\ A. 
Brudny{\u \i}, N. Ya.\ Kruglyak, and E. I. Pustylnik for valuable 
discussions and our warmest appreciation to M. Cwikel for his interest 
in this research.

\enddocument